\documentclass[11pt, twoside, eqno]{article}
\usepackage{latexsym}
\usepackage{amssymb}
\usepackage{amsfonts}
\begin{document}
\begin{center}

\noindent {\bf \Large Comment on "Generalized ideal elements in
$le$-$\Gamma$-semigroups" }\bigskip

\medskip

{\bf Niovi Kehayopulu, Michael Tsingelis}\bigskip
\end{center}
\bigskip

\noindent Concerning the paper in the title by K. Hila and E. Pisha 
in Commun. Korean Math. Soc. Volume 26, Issue 3 (2011), 373--384 [1], 
we give our results and make the main corrections. The Theorem 2.4, 
Theorem 2.5, Theorem 2.6, Theorem 3.4, Theorem 3.9 of the paper are 
based on Lemma 2.3, but Lemma 2.3
is wrong. As Theorem 3.4 is wrong, Theorem 3.5 is wrong as well. The
$(m,0)$-ideal elements and the $(0,n)$-ideal elements mentioned in 
Theorem 9.6 have not defined in the paper (look at Definition 2.1). 
Lemma 3.7 is also wrong as the expression $a^n$ is used in it. In 
Lemma
2.3(3), the authors use the Definition 2.2 which clearly is not
true: Let $M$ be a $\vee e$-$\Gamma$-semigroup and $m,n\in Z^{+}$.
According to the paper, $M$ is called $(m,n)$-regular if for all
$a\in M$ and all $\lambda, \mu\in\Gamma$ there exist
$\gamma_1,\gamma_2, ..... ,\gamma_{m-1},\rho_1,\rho_2, .....
,\rho_{n-1}\in\Gamma$ such that$$a\le (a\gamma_1 a\gamma_2 a .....
\gamma_{m-1}a)\lambda e\mu (a\rho_1,a\rho_2, .....
,\rho_{n-1}a).$$Suppose $m=1$ (or $n=1$). What is the
$\gamma_{m-1}$ (that is, the $\gamma_0$) in the expression
$a\gamma_1 a\gamma_2 a ..... \gamma_{m-1}a$ ? What is the
$\rho_{m-1}=\rho_0$ in the expression $a\rho_1 a\rho_2 a .....
\rho_{m-1}a$ ? The Definition 2.2 of the paper should be corrected
and then to check if its corrected form coincides with the
definition of regular (that is, $(1,1)$-regular) $\vee
e$-$\Gamma$-semigroups. Lemma 2.3(3) is based on the definition of
$(m,n)$-ideal elements as well given in Definition 2.1 of the
paper, which is also wrong. According to Definition 2.1, an
element $a$ of a $\vee e$-$\Gamma$-semigroup is an $(m,n)$-ideal
element ($m,n\in Z^{+}$) if there exist $\gamma_1,\gamma_2, .....
,\gamma_{m-1},\rho_1,\rho_2, ..... ,\rho_{n-1}\in\Gamma$ such
that$$(a\gamma_1 a\gamma_2 a ..... \gamma_{m-1}a)\lambda e\mu
(a\rho_1,a\rho_2, ..... ,\rho_{n-1}a)\le a$$ for all
$\lambda,\mu\in\Gamma$. The same question arises: Suppose $m=1$
(or $n=1$). What is the $\gamma_{m-1}$ (that is, the $\gamma_0$)
in the expression $a\gamma_1 a\gamma_2 a ..... \gamma_{m-1}a$ ?
What is the $\rho_{m-1}=\rho_0$ in the expression $a\rho_1 a\rho_2
a ..... \rho_{m-1}a$ ? The Definition 2.1 should be corrected and
then one has to check if its corrected form coincides with the
definition of a bi-ideal element (that is, $(1,1)$-ideal element).
As far as the definition of a quasi-ideal element and the
definition of bi-ideal element is concerned the authors gave the
following definitions: An element $a$ of a $\vee
e$-$\Gamma$-semigroup is called a quasi-ideal element if, for  for 
every $\lambda\in\Gamma$, $a\lambda
e\wedge e\lambda a$ exists and $a\lambda e\wedge e\lambda a\le a$.
The element $a$ is called a bi-ideal element if $a\lambda e\mu
a\le a$ for all $\lambda,\mu\in\Gamma$. But the quasi-ideal
elements should be bi-ideal elements as well. There is no such a
proof in the paper, and it does not seem to be true. One can
construct an example using tables which shows that this is not
true.

For each of the results of the papers in [2],[3], the authors
tried to get its analogous in case of a $\vee
e$-$\Gamma$-semigroup just, casually, putting $\alpha$, $\beta$
(the elements of $\Gamma$) in some places. For shortly, they wrote
$a^m$ as the element $a\gamma_1 a\gamma_2 a ..... \gamma_{m-1}a$
for some $\gamma_1,\gamma_2, ..... , \gamma_{m-1}\in\Gamma$ $(m\in
Z^+$) which leads to the mistakes throughout the paper. Look, for
example, at Lemma 2.3(1),(2). Besides, in the proof of Lemma
2.3(1) the authors wrote, $a\gamma e\vee a^m \lambda e\mu a^k
\gamma e=a\gamma e$, $a\rho a\gamma e\vee a^m \lambda e\mu
a^{k+1}\gamma e=a^2\gamma e$ which certainly is not true. Except
of the fact that $a^m$, $a^n$ has been used in Lemma 2.3(3), this 
part of the
lemma has an additional mistake. According to Lemma 2.3(3),
$<a>_{(m,n)}=a\vee a^m\lambda e\mu a^n$ for every
$\lambda,\mu\in\Gamma$. $<a>_{(m,n)}$ is uniquely defined, while
they consider it equal to $a\vee a^m\lambda e\mu a^n$ for every
$\lambda,\mu\in\Gamma$. If this is the case, the authors should prove 
that for every $\lambda,\mu\in\Gamma$, $a\vee a^m\lambda e\mu a^n$ is 
uniquely defined. Is it possible ? That is, if
$\lambda,\mu\in\Gamma$ and $\gamma,\delta\in\Gamma$ then is $a\vee
a^m\lambda e\mu a^n=a\vee a^m\gamma e\delta a^n$ ? Let us get
$m=2$, $n=2$, for example. According to the paper by Hila and
Pisha, an element $a$ of $M$ is called a $(2,2)$-ideal element if
there exist $\gamma,\delta\in M$ such that $(a\gamma a)\xi e\zeta
(a\delta a)\le a$ for all $\xi,\zeta\in\Gamma$. Then they write
$a^2\xi e\zeta a^2\le a$, which actually means that $a\gamma
a=a\delta a$. In that case they should prove that $a\gamma
a=a\delta a$. Is it so, and why? So they cannot write
$<a>_{(2,2)}=a\vee a^2\xi e\zeta a^2$ for all
$\xi,\zeta\in\Gamma$. Shortly, Lemma 2.3 is without any sense, and
so is the rest of the paper. Finally, the definition of a $\vee
e$-$\Gamma$-semigroup given in Definition 1.7 is also not correct.
The authors say: Let "$M$ be a semilattice under $\vee$ ..." which
means that there exists an order relation $\le$ on $M$ according
to which $M$ is a semilattice, that is, for any two elements
$a,b\in M$ there exists an element $t\in M$ (denoted by $a\vee b$
and called the supremum of $a$ and $b$) such that $t\ge a$, $t\ge
b$ and if $h\in M$ such that $h\ge a$ and $h\ge b$, then $t\le h$
($y\ge x$ means $x\le y$ that is $(x,y)\in \le$). Then, they say

"The usual order relation $"\le"$ on $M$ is defined in the
following way$$a\le b \;\Longleftrightarrow\; a\vee b=b."$$ and
they add that $a\le b$ implies $a\gamma c\le b\gamma c$ and
$c\gamma a\le c\gamma b$ for all $c\in M$ and all
$\gamma\in\Gamma$. This has no sense, it is wrong, because the
order defines the semilattice. Besides, they should mention that
$a\le b$ implies $a\gamma c\le b\gamma c$ and $c\gamma a\le
c\gamma b$ for all $c\in M$ and all $\gamma\in\Gamma$ immediately
after the Definition 1.7 in line 20 of page 375 and not on lines
22-23 as they did. Its proof is as follows: Let $M$ be a $\vee
e$-$\Gamma$-semigroup, $a\le b$, $\gamma\in\Gamma$ and $c\in M$.
Since $a\le b$, we have $a\vee b=b$. Then, by Definition 1.1(3) of
the paper, we have $(a\vee b)\gamma c=b\gamma c$. Since $M$ is a
$\vee e$-$\Gamma$-semigroup, $(a\vee b)\gamma c=a\gamma c\vee
b\gamma c$. Then $a\gamma c\le a\gamma c\vee b\gamma c=b\gamma c$,
and $a\gamma c\le b\gamma c$. Similarly $a\le b$ implies $c\gamma
a\le c\gamma b$ for all $c\in M$ and all $\gamma\in\Gamma$. This
means that every $\vee e$-$\Gamma$-semigroup is a
$poe$-$\Gamma$-semigroup. In addition, instead of "for all
$a,b,c\in M$" written in the paper, is much better to write for
all $c\in M$ since $a\le b$ that is $(a,b)\in\le (\subseteq
M\times M)$ implies $a,b\in M$.

The authors tried to extend the results given in [2],[3] from
$\vee e$-semigroups to ordered $\Gamma$-semigroups, but there is 
nothing correct in this paper as the expression $a^n$ has been used 
throughout the paper.\smallskip

In the following we correct the results given by the authors in
Lemma 2.3, Theorem 2.4 and Theorem 2.5 which correspond to the
result in [2] (Lemma 1, Theorem 1, Theorem 2 in [2]). Based on our 
results, the authors might correct the
rest of their paper which corresponds to the results given in [3].

The first two properties of Lemma 2.3 might be corrected as
follows:\medskip

\noindent{\bf Definition 1.} Let $M$ be a $\Gamma$-semigroup,
$a\in M$, $\gamma\in\Gamma$ and $m\in N$ ($N=\{1,2, ..... , n\}$
is the set of natural numbers). Then

(1) if $m=1$, we define $a_{\gamma}^1:=a$

(2) if $m\ge 2$, we define $a_\gamma^m:=a\gamma a\gamma
\,.....\,\gamma a$

\hspace{2cm} ($m-1$-times the $\gamma$, $m$-times the
$a$).\medskip

\noindent{\bf Remark 2.} If $a\in M$, $\gamma\in\Gamma$ and
$m,n\in N$, then we have $$a_\gamma^m \gamma
a_\gamma^n=a_\gamma^{m+n}.$$The first two properties of Lemma 1 in
[1] can be formulated as follows:\medskip

\noindent{\bf Lemma 3.} {\it Let M be a $\vee
e$-$\Gamma$-semigroup, $a,b\in M$, $\gamma\in\Gamma$ and $m,n\in N$.
Then we have

$(1)$ $(a\vee a_\gamma^m \gamma b)_\gamma^m\gamma e=a_\gamma^m \gamma 
e$.

$(2)$ $e\gamma (a\vee b\gamma a_\gamma^n)_\gamma^n=e\gamma
a_\gamma^n$.}\medskip

\noindent{\bf Proof.} (1) For $m=1$, condition (1) is
satisfied. Indeed: Let $a,b\in M$ and $\gamma\in\Gamma$. Then we
have\begin{eqnarray*}(a\vee a_\gamma^m \gamma b)_\gamma^m\gamma
e&=&(a\vee a_\gamma^1 \gamma b)_\gamma^1\gamma e=(a\vee a\gamma 
b)\gamma
e\\&=&a\gamma e\vee a\gamma b\gamma e.\end{eqnarray*}Since
$b\gamma e\le e$, we have $a\gamma b\gamma e\le a\gamma e$, then
$a\gamma e\vee a\gamma b\gamma e=a\gamma e$, and$$(a\vee
a_\gamma^m \gamma b)_\gamma^m\gamma e=a\gamma e=a_\gamma^1\gamma
e=a_\gamma^m\gamma e.$$Suppose condition (1) is satisfied for
$m=k\ge 1$. That is, suppose that for every $c,d\in M$ and every
$\delta\in\Gamma$, we have $$(c\vee c_\delta^k \delta d)_\delta^k
\delta e=c_\delta^k \delta e.$$Then it is satisfied
for $m=k+1$ as well. That is, for every $a,b\in M$ and every
$\gamma\in\Gamma$, we have$$(a\vee a_\gamma^{k+1}\gamma
b)_\gamma^{k+1}\gamma e=a_\gamma^{k+1}\gamma e.$$Indeed: Let
$a,b\in M$ and $\gamma\in\Gamma$. Then we
have\begin{eqnarray*}(a\vee a_\gamma^{k+1}\gamma
b)_\gamma^{k+1}\gamma e&=&\Big ( (a\vee a_\gamma^{k+1}\gamma
b)_\gamma^1 \gamma(a\vee a_\gamma^{k+1}\gamma b)_\gamma^k)\Big )
\gamma e\;\;\mbox {(by Remark 2)}\\&=&(a\vee a_\gamma^{k+1}\gamma
b)_\gamma^1 \gamma \Big ( (a\vee a_\gamma^{k+1}\gamma b)_\gamma ^k
\gamma e\Big )\mbox { (since } M \mbox { is
associative)}\\&=&(a\vee a_\gamma^{k+1}\gamma b)\gamma\Big (
(a\vee a_\gamma^{k+1}\gamma b)_\gamma^k \gamma e\Big )\mbox { (by
Definition 1(1))}\\&=&(a\vee a_\gamma^{k+1}\gamma b)\gamma\bigg (
\Big (a\vee (a_\gamma^k\gamma a_\gamma^1)\gamma b\Big )_\gamma^k
\gamma e\bigg )\mbox { (by Remark 2)}\\&=&(a\vee
a_\gamma^{k+1}\gamma b)\gamma\bigg ( \Big (a\vee (a_\gamma^k\gamma
a)\gamma b\Big )_\gamma^k \gamma e\bigg )\mbox { (by Definition
1(1)}\\&=&(a\vee a_\gamma^{k+1}\gamma b)\gamma \bigg (\Big (a\vee
a_\gamma^k\gamma (a\gamma b)\Big )_\gamma^k \gamma e\bigg )\mbox {
(since } M \mbox { is associative)}\\&=&(a\vee
a_\gamma^{k+1}\gamma b)\gamma (a_\gamma^k\gamma e) \mbox { (by the
assumption})\\&=&(a\gamma a_\gamma^k \gamma e)\vee
(a_\gamma^{k+1}\gamma b\gamma a_\gamma^k \gamma e)\mbox { (since }
M \mbox { is a $\vee e$-$\Gamma$-semigroup)}
\\&=&(a_\gamma^{k+1}\gamma e)\vee (a_\gamma^{k+1}\gamma b\gamma
a_\gamma^k \gamma e)\mbox { (by Remark 2)}.\end{eqnarray*}Since
$b\gamma a_\gamma^k \gamma e\le e$, we have $a_\gamma^{k+1}\gamma 
b\gamma a_\gamma^k \gamma e\le a_\gamma^{k+1}\gamma e$, and so
$$a_\gamma^{k+1}\gamma e\vee a_\gamma^{k+1}\gamma b\gamma
a_\gamma^k\gamma e=a_\gamma^{k+1}\gamma e.$$Thus we have $(a\vee
a_\gamma^{k+1}\gamma b)_\gamma^{k+1}\gamma e=a_\gamma^{k+1}\gamma
e$. Condition (2) can be proved in a similar
way.$\hfill\Box$\medskip

\noindent {\bf Corollary 4.} {\it Let M be a $\vee
e$-$\Gamma$-semigroup, $a,b\in M$, $\beta,\gamma\in\Gamma$ and
$m,n\in N$. Then we have

$(1)$ $(a\vee a_\beta^m \beta e\gamma a_\gamma^n)_\beta^m\beta
e=a_\beta^m \beta e$.

$(2)$ $e\gamma (a\vee a_\beta^m\beta e\gamma
a_\gamma^n)_\gamma^n=e\gamma a_\gamma^n$.}\medskip

\noindent{\bf Proof.} (1) Since $M$ is associative, we have
$(a\vee a_\beta^m \beta e\gamma a_\gamma^n)_\beta^m\beta e=\Big
(a\vee a_\beta^m \beta (e\gamma a_\gamma^n )\Big )_\beta^m\beta
e.$ We put $e\gamma a_\gamma^n=c$ and, by Lemma 3(1), we have
$$(a\vee a_\beta^m \beta e\gamma a_\gamma^n)_\beta^m\beta e=(a\vee
a_\beta^m\beta c)_\beta^m\beta e=a_\beta^m \beta e.$$(2) We have
$e\gamma (a\vee a_\beta^m\beta e\gamma
a_\gamma^n)_\gamma^n=e\gamma\Big (a\vee (a_\beta^m\beta e)\gamma
a_\gamma^n\Big )_\gamma^n$. We set $a_\beta^m\beta e=d$ and, by
Lemma 3(2), we get $e\gamma (a\vee a_\beta^m\beta e\gamma
a_\gamma^n)_\gamma^n=e\gamma (a\vee d\gamma
a_\gamma^n)_\gamma^n=e\gamma a_\gamma^n$.$\hfill\Box$\smallskip

\noindent As far as the third property of Lemma 1 in [1] is
concerned, we first have to introduce the following
definition:\medskip

\noindent{\bf Definition 5.} Let $M$ be a $\vee
e$-$\Gamma$-semigroup, $m,n\in N$ and $\beta,\gamma\in\Gamma$. An
element $a$ of $M$ is called an {\it $(m,n,\beta,\gamma)$-ideal
element} if$$a_\beta^m\beta e\gamma a_\gamma^n\le a.$$For
$\beta=\gamma$, the element $a$ is called an {\it $(m,n,\beta)$-ideal 
element}.
An element $a$ of $M$ is called an {\it $(m,0,\beta)$-ideal
element} if$$a_\beta^m\beta e\le a.$$It is called an {\it
$(0,m,\beta)$-ideal element} if $e\beta
a_\beta^m\le a$. For $a\in M$ denote by
$<a>_{(m,n,\beta,\gamma)}$ the $(m,n,\beta,\gamma)$-ideal element
of $M$ generated by $a$; denote by $<a>_{(m,0,\beta)}$ the
$(m,0,\beta)$-ideal element of $M$ generated by $a$ and
$<a>_{(0,m,\beta)}$ the $(0,m,\beta)$-ideal element of $M$
generated by $a$. We denote by $I_{(m,n,\beta,\gamma)}$ the set of
$(m,n,\beta,\gamma)$-ideal elements of $M$ and by
$I_{(m,0,\beta)}$ (resp. $I_{(0,m,\beta)})$ the set of
$(m,0,\beta)$ (resp. $(0,m,\beta)$)-ideal elements of $M$.
\smallskip

\noindent The Lemma 2.3(3) in [1] should be formulated as
follows:\medskip

\noindent{\bf Lemma 6.} {\it Let M be a $\vee
e$-$\Gamma$-semigroup, $a\in M$, $m,n\in M$ and
$\beta,\gamma\in\Gamma$. Then we
have$$<a>_{(m,n,\beta,\gamma)}=a\vee a_\beta^m\beta e\gamma
a_\gamma^n.$$}{\bf Proof.} The element $a\vee a_\beta^m\beta
e\gamma a_\gamma^n$ is an $(m,n,\beta,\gamma)$-ideal element of
$M$. That is,$$(a\vee a_\beta^m\beta e\gamma
a_\gamma^n)_\beta^m\beta e\gamma (a\vee a_\beta^m\beta e\gamma
a_\gamma^n)_\gamma^n\le a\vee a_\beta^m\beta e\gamma
a_\gamma^n.$$Indeed,\begin{eqnarray*}(a\vee a_\beta^m\beta e\gamma
a_\gamma^n)_\beta^m\beta e\gamma (a\vee a_\beta^m\beta e\gamma
a_\gamma^n)_\gamma^n&=&\Big ((a\vee a_\beta^m\beta e\gamma
a_\gamma^n)_\beta^m \beta e\Big )\gamma (a\vee a_\beta^m\beta
e\gamma a_\gamma^n)_\gamma^n\\&=&(a_\beta^m\beta e)\gamma (a\vee
a_\beta^m\beta e\gamma a_\gamma^n)_\gamma^n\mbox { (by Corollary
4(1))}\\&=&a_\beta^m\beta \Big (e\gamma (a\vee a_\beta^m\beta
e\gamma a_\gamma^n)_\gamma^n\Big )\\&=&a_\beta^m\beta e\gamma
a_\gamma^n\mbox { (by Corollary 4(2))}\\&\le&a\vee a_\beta^m\beta
e\gamma a_\gamma^n.\end{eqnarray*}Clearly, $a\le a\vee
a_\beta^m\beta e\gamma a_\gamma^n$. Let now $t$ be an
$(m,n,\beta,\gamma)$-ideal element of $M$ such that $t\ge a$. Then
$a\vee a_\beta^m\beta e\gamma a_\gamma^n\le t$. Indeed: Since
$t\ge a$, we have $a_\beta^m\le t_\beta^m$ and $a_\gamma^n\le
t_\gamma^n$. Then we have $a\vee a_\beta^m\beta e\gamma
a_\gamma^n\le a\vee t_\beta^m\beta e\gamma t_\gamma^n\le
t$.$\hfill\Box$\medskip

\noindent{\bf Definition 7.} Let $M$ be a $\vee
e$-$\Gamma$-semigroup, $m,n\in M$ and $\beta,\gamma\in\Gamma$. An
element $a$ of $M$ is called {\it $(m,n,\beta,\gamma)$-regular} if
$a\le a_\beta^m\beta e\gamma a_\gamma^n$. $M$ is called {\it
$(m,n,\beta,\gamma)$-regular} if every element of $M$ is so.
\medskip

\noindent{\bf Theorem 8.} {\it Let $M$ be a $\vee
e$-$\Gamma$-semigroup, $m,n\in M$ and $\beta,\gamma\in\Gamma$. The
following are equivalent:

$(1)$ M is $(m,n,\beta,\gamma)$-regular.

$(2)$ $a_\beta^m\beta e\gamma a_\gamma^n=a$ for every $a\in
I_{(m,n,\beta,\gamma)}$. }\medskip

\noindent{\bf Proof.} $(1)\Longrightarrow (2)$. Let $a\in
I_{(m,n,\beta,\gamma)}$. Then $a_\beta^m\beta e\gamma
a_\gamma^n\le a$. Since $a\in M$, by hypothesis, we have $a\le
a_\beta^m\beta e\gamma a_\gamma^n$. Thus we have $a_\beta^m\beta
e\gamma a_\gamma^n=a$.\smallskip

\noindent$(2)\Longrightarrow (1)$. Let $a\in M$. Then $a\le
a_\beta^m\beta e\gamma a_\gamma^n$. Indeed: Since
$<a>_{(m,n,\beta,\gamma)}$ is an $(m,n,\beta,\gamma)$-ideal
element of $M$, by hypothesis, we
have$$\Big(<a>_{(m,n,\beta,\gamma)}\Big)_\beta^m\beta e\gamma
\Big(<a>_{(m,n,\beta,\gamma)}\Big)_\gamma^n=<a>_{(m,n,\beta,\gamma)}.$$
Then, by Lemma 6, we get\begin{eqnarray*}a&\le& a\vee
a_\beta^m\beta e\gamma a_\gamma^n=(a\vee a_\beta^m\beta e\gamma
a_\gamma^n)_\beta^m\beta e\gamma (a\vee a_\beta^m\beta e\gamma
a_\gamma^n)_\gamma^n\\&=&\Big((a\vee a_\beta^m\beta e\gamma
a_\gamma^n)_\beta^m\beta e\Big)\gamma (a\vee a_\beta^m\beta
e\gamma a_\gamma^n)_\gamma^n\mbox { (since } M \mbox { is
associative)}\\&=&(a_\beta^m\beta e)\gamma (a\vee a_\beta^m\beta
e\gamma a_\gamma^n)_\gamma^n\mbox { (by Corollary 4(1))
}\\&=&a_\beta^m\beta \Big(e\gamma (a\vee a_\beta^m\beta e\gamma
a_\gamma^n)_\gamma^n\Big)\mbox { (since } M \mbox { is
associative)}\\&=&a_\beta^m\beta e\gamma a_\gamma^n\mbox { (by
Corollary 4(2))},\end{eqnarray*}and the proof is
complete.$\hfill\Box$\medskip

\noindent{\bf Definition 9.} An element $a$ of a
$po$-$\Gamma$-semigroup $M$ is called {\it $\beta$-subidempotent} if
$a\beta a\le a$. $M$ is called {\it $\beta$-subidempotent} if every
element of $M$ is so.\medskip

\noindent{\bf Lemma 10.} {\it Let M be a $\vee
e$-$\Gamma$-semigroup, $a\in M$, $m\in M$ and $\beta\in\Gamma$.
Then we have

$(1)$ $<a>_{(m,0,\beta)}=a\vee a_\beta^m\beta e$.

$(2)$ $<a>_{(0,m,\beta)}=a\vee e\beta a_\beta^m$.} \medskip

\noindent {\bf Proof.} (1) By Lemma 3(1), we have $(a\vee 
a_\beta^m\beta e)_\beta^m\beta
e=a_\beta^m\beta e\le a\vee a_\beta^m\beta e$, so $a\vee
a_\beta^m\beta e$ is an $(m,0,\beta)$-ideal element of $M$
containing $a$. If now $t$ is an $(m,0,\beta)$-ideal element of
$M$ such that $a\le t$, then we have $a\vee a_\beta^m\beta e\le
t\vee t_\beta^m\beta e=t$. The proof of (2) is similar.
$\hfill\Box$\smallskip

\noindent By Lemma 3, the following lemma holds:\medskip

\noindent {\bf Lemma 11.} {\it Let M be a $\vee
e$-$\Gamma$-semigroup, $a,b\in M$, $\beta,\gamma\in\Gamma$ and
$m,n\in N$. Then we have

$(1)$ $(a\vee a_\beta^m \beta e)_\beta^m\beta e=a_\beta^m \beta
e$.

$(2)$ $e\gamma (a\vee e\gamma a_\gamma^n)_\gamma^n=e\gamma
a_\gamma^n$.}\medskip

\noindent{\bf Theorem 12.} {\it Let M be an $\ell
e$-$\Gamma$-semigroup, $m,n\in M$ and $\beta,\gamma\in\Gamma$.
Suppose M is $\beta$-subidempotent. The following are equivalent:

$(1)$ M is $(m,n,\beta,\gamma)$-regular.

$(2)$ $a\wedge b=a_\beta^m\beta b\wedge a\gamma b_\gamma^n$ for
every $a\in I_{(m,0,\beta)}$ and every $b\in
I_{(0,n,\gamma)}$.}\medskip

\noindent{\bf Proof.} $(1)\Longrightarrow (2)$. Let $a\in
I_{(m,0,\beta)}$ and $b\in I_{(0,n,\gamma)}$. Since $a\in
I_{(m,0,\beta)}$, we have $a_\beta^m\beta e\le a$, then
$a_\beta^m\beta b\le a_\beta^m\beta e\le a$. Since $b\in
I_{(0,n,\gamma)}$, we have $e\gamma b_\gamma^n\le b$, then
$a\gamma b_\gamma^n\le e\gamma b_\gamma^n\le b$. Thus we have
$a_\beta^m\beta b\wedge a\gamma b_\gamma^n\le a\wedge b$. Since
$a\wedge b\in M$ and $M$ is $(m,n,\beta,\gamma)$-regular, we have
$a\wedge b\le (a\wedge b)_\beta^m\beta e\gamma (a\wedge
b)_\gamma^n$. Since $a\wedge b\le a$, we have $(a\wedge
b)_\beta^m\le a_\beta^m$. Since $a\wedge b\le b$, we have
$(a\wedge b)_\gamma^n\le b_\gamma^n$. Then$$a\wedge b\le (a\wedge
b)_\beta^m\beta e\gamma (a\wedge b)_\gamma^n\le a_\beta^m\beta
e\gamma b_\gamma^n.$$Since $a_\beta^m\beta (e\gamma b_\gamma^n)\le
a_\beta^m\beta b$ and $(a_\beta^m\beta e)\gamma b_\gamma^n\le
a\gamma b_\gamma^n$, we have$$a\wedge b\le a_\beta^m\beta e\gamma
b_\gamma^n\le a_\beta^m\beta b\wedge a\gamma b_\gamma^n,$$and
property (2) is satisfied.\smallskip

\noindent$(2)\Longrightarrow (1)$. We remark first that $M$ has the 
following property:

$$c\wedge d\le c\beta d \mbox { for every } c\in I_{(m,0,\beta)}
\mbox { and every } d\in I_{(0,n,\gamma)} \;..........\;
(*)$$Indeed: Let $c\in I_{(m,0,\beta)}$ and $d\in
I_{(0,n,\gamma)}$. By (2), we have$$c\wedge d\le c_\beta^m\beta
d\wedge c\gamma d_\gamma^n\le  c_\beta^m\beta d.$$Since $M$ is
$\beta$-subidempotent, we have $c_\beta^m\le c$, so $c\wedge d\le
c\beta d$, and $(*)$ holds.\\Let now $a\in M$. Then $a\le
a_\beta^m\beta e\gamma a_\gamma^n$. Indeed: $$a\le
<a>_{(m,0,\beta)}\wedge <a>_{(0,n,\gamma)}\le
<a>_{(m,0,\beta)}\beta <a>_{(0,n,\gamma)} \mbox { (by } 
(*)).$$Moreover,\begin{eqnarray*}<a>_{(m,0,\beta)}&=&<a>_{(m,0,\beta)}\wedge 
e\\&=&\Big (<a>_{(m,0,\beta)}\Big )_\beta^m\beta e\wedge 
<a>_{(m,0,\beta)}\gamma e_\gamma^n \mbox { (by (2)})\\&\le& \Big 
(<a>_{(m,0,\beta)}\Big )_\beta^m\beta e\\&=&(a\vee a_\beta^m\beta 
e)_\beta^m\beta e \mbox { (by Lemma 10)}\\&=&a_\beta^m\beta e \mbox { 
(by Lemma 11).}
\end{eqnarray*}On the other hand, by Lemma 10, we have 
$<a>_{(m,0,\beta)}\ge a_\beta^m\beta e$. Thus we have 
$<a>_{(m,0,\beta)}=a_\beta^m\beta e$. Similarly we have 
$<a>_{(0,n,\gamma)}=e\gamma a_\gamma^n$. Therefore, we obtain$$a\le 
(a_\beta^m\beta e)\beta (e\gamma a_\gamma^n)=a_\beta^m\beta
(e\beta e)\gamma a_\gamma^n\le
a_\beta^m\beta e\gamma a_\gamma^n,$$and $M$ is 
$(m,n,\beta,\gamma)$-regular.$\hfill\Box$
 \bigskip

\noindent This paper has been submitted in Commun. Korean Math. Soc. 
on October 31, 2013 at 17:15 (the date in Greece).

\end{document}